\documentclass[10pt,a4paper]{amsart}

\usepackage{latexsym,amssymb,epsfig,array,color}
\usepackage{mathrsfs,mathabx}
\usepackage{tikz, tikz-cd}
\usepackage{graphicx,epstopdf}
\usepackage{colortbl}
\usepackage[all]{xy}
\usepackage{tipa}
\usepackage{wrapfig}

\theoremstyle{plain}
\newtheorem{theorem}{Theorem}
\newtheorem{prop}{Proposition}
\newtheorem{corollary}{Corollary}
\newtheorem{corollary_th}{Corollary}[theorem]
\newtheorem{lemma}{Lemma}

\theoremstyle{remark}
\newtheorem{remark}{Remark}
\newtheorem{definition}{Definition}

\begin{document}

\author{Florian L. Deloup}
\title{On Heisenberg groups}

\begin{abstract}
It is known that an abelian group $A$ and a $2$-cocycle $c:A \times A \to C$ yield a group ${\mathscr{H}}(A,C,c)$ which we call a Heisenberg group. This group, a central extension of $A$, is the archetype of a class~$2$ nilpotent group. In this note, we prove that under mild conditions, any class~$2$ nilpotent group $G$ is equivalent as an extension of $G/[G,G]$ to a Heisenberg group ${\mathscr{H}}(G/[G,G], [G,G], c')$ whose $2$-cocycle $c'$ is bimultiplicative. 
\end{abstract}

\maketitle
\tableofcontents


\section{Introduction}

The Heisenberg group is a famous noncommutative group, possibly one of the very first examples one encounters in linear algebra along with transformation groups. Various incarnations of it appear in many fields, including harmonic analysis, complex analysis, representation theory and quantum topology. Its simplest incarnation $H({\mathbb{R}})$ is the group that consists of upper triangular $3 \times 3$  matrices with $1$'s on the diagonal with usual matrix product. It is a noncommutative nilpotent subgroup of the group ${\rm{GL}}_{3}(\mathbb{R})$ of invertible $3 \times 3$ matrices. 

Here we consider a certain generalized version of the Heisenberg group $H_{\beta}(A) = C \times A$ where $\beta:A \times A \to C$ is a bilinear pairing on an Abelian group $A$ with group law:
$$ (t,x) \cdot (t',y) = (t+t'+\beta(x,y), x+y), \ \  \ t,t' \in C, \ x, y \in A.$$
This Heisenberg group $H_{\beta}(A)$ plays a basic and important r\^ole in Abelian topological quantum field theories \cite{Del}. The main result of this note is that, under mild conditions on the commutator subgroup $[G,G]$ and the abelianized group $G/[G,G]$, any nilpotent group $G$ of nilpotency class $2$ is equivalent, as an extension, to some Heisenberg group $H_{\beta}(A)$. Furthermore, such extensions are classified by the symplectic pairing $\omega_{\beta}:G/[G,G] \times G/[G,G] \to [G,G]$ with $\omega_{\beta}(x,y) = \beta(x,y) - \beta(y,x)$, $x, y \in G/[G,G]$. 
Foundational consequences for Abelian Topological Quantum Field Theories are discussed elsewhere \cite{Del}.


\section{Initial set-up and main theorem}

\begin{definition}
Let $A$ and $C$ be Abelian groups (written multiplicatively) and let $\beta:A \times A \to C$ be a bimultiplicative pairing. 
The {\emph{Heisenberg group}} ${\mathscr{H}}_{\beta}(A)$ is the extension of $A$ defined as the set $C \times A$ endowed
with the multiplication rule
$$ (t,x) \cdot (t',y) = (tt'\beta(x,y), xy).$$
If $\beta$ is understood, we suppress the subscript and write simply ${\mathscr{H}}(A)$.
\end{definition}

Associativity follows from associativity in $A$ and bimultiplicativity of $\beta$; the pair $(1,1)$ of neutral elements in $C$ and $A$ respectively is the neutral element of ${\mathscr{H}}(A)$; the inverse of $(t,x)$ is $(t^{-1}\beta(x,x),x^{-1})$ for any $x \in A$ and $t \in {C}$.

It follows from the definition that ${\mathscr{H}}(A)$ lies in the exact sequence of groups
\begin{equation}
\xymatrix{ 1 \ar[r] & {C} \ar[r] & {\mathscr{H}}(A)
\ar[r] & A \ar[r] & 1}
\label{eq:exact_sequence_heisenberg}
\end{equation}
where $C \to {\mathscr{H}}(A), t \mapsto (t,1)$ is the natural inclusion and ${\mathscr{H}}(A) \to A, (t,x) \mapsto x$ is the projection onto the second factor.

\begin{definition}
Given a bimultiplicative pairing $\beta:A \times A \to C$, the pairing $\omega_{\beta}:A \times A \to C$ defined by
\begin{equation}
\omega_{\beta}(x,y) = \beta(x,y)\cdot \beta(y,x)^{-1}, \ \ \ x,y \in A,
\label{eq:omega}
\end{equation}
is the symplectic pairing associated to $\beta$. When $\beta$ is self-understood, we suppress the subscript and write simply $\omega$.
We say that $\omega$ is nondegenerate (resp. regular) if the adjoint map $\hat{\omega}: A \to {\rm{Hom}}(A,C), \ x \mapsto \omega(x,-)$ is injective (resp. bijective).
\end{definition}

\begin{lemma}
For any bimultiplicative pairing $\beta$, $\omega_{\beta}$ is alternating (hence antisymmetric).
\end{lemma}

\begin{proof}
Obvious from the definition (\ref{eq:omega}).
\end{proof}

\begin{lemma} \label{lem:Heisenberg-observations}
Let $X = (t,x), Y = (t',y) \in \mathscr{H}(A)$. Then
\begin{equation}
 [X,Y] = (\omega(x,y), 1).\label{eq:Heisenberg-commutator}
 \end{equation}
\end{lemma}

\begin{proof}
The equality is a direct computation:
\begin{align*}
[X,Y] & = XYX^{-1}Y^{-1} \\
& = (t,x)(t',y)(t^{-1}\beta(x,x),x^{-1})(t'^{-1}\beta(y,y),y^{-1}) \\
& = (tt'\beta(x,y), xy)(t^{-1}t'^{-1}\beta(x,x)\beta(y,y)\beta(x^{-1},y^{-1}), x^{-1}y^{-1}) \\
& = (\beta(x,y)\beta(x,x)\beta(y,y)\beta(x,y)\beta(xy,xy)^{-1},1) \\
& = (\beta(x,y)\beta(y,x)^{-1},1) \\
& = (\omega(x,y),1).
\end{align*}
\end{proof}

\begin{lemma}
The commutator subgroup $[{\mathscr{H}}(A), {\mathscr{H}}(A)]$ is $C_{\omega} \times 1$ where $C_{\omega}$ is the subgroup of $C$ generated by the image of $\omega$.
The center $Z({\mathscr{H}}(A))$ of ${\mathscr{H}}(A)$ is $C \times {\rm{Ker}}\ \hat{\omega}$. In particular, 
$$ [{\mathscr{H}}(A), {\mathscr{H}}(A)] \subseteq C \times 1 \subseteq Z({\mathscr{H}}(A)).$$
\end{lemma}

\begin{proof}
Both statements follow from Lemma \ref{lem:Heisenberg-observations}.
\end{proof}

\begin{corollary}
The symplectic pairing $\omega$ is nondegenerate if and only if  $Z({\mathscr{H}}(A)) = C \times 1$. If furthermore $C$ is generated by the image of $\omega$, then
$$ [{\mathscr{H}}(A), {\mathscr{H}}(A)] = C \times 1= Z({\mathscr{H}}(A)).$$
\end{corollary}

\begin{corollary} \label{cor:heisenberg_is_2-nilpotent}
The Heisenberg group ${\mathscr{H}}(A)$ is nilpotent of nilpotency class at most two.
\end{corollary}

The goal of this note is to prove a strong converse of Corollary \ref{cor:heisenberg_is_2-nilpotent}.

The following result should be well-known.

\begin{prop} \label{prop:nilpotent2}
Let $G$ be a group, $[G,G] \triangleleft G$ its commutator subgroup and $Z(G) \triangleleft G$ its center. The following assertions are equivalent:
\begin{enumerate}
\item[(1)] The nilpotency class of $G$ is $2$;
\item[(2)] $[G,G] \subseteq Z(G)$;
\item[(3)] The commutator map $G \times G \to [G,G], (g,h) \mapsto [g,h] = ghg^{-1}h^{-1}$ descends to an alternating nondegenerate bimultiplicative pairing
$$ \varpi: G/Z(G) \times G/Z(G) \to [G,G], \ ([g], [h]) \mapsto [g,h].$$
\end{enumerate}
\end{prop}

\begin{proof}
For an arbitrary group $G$, the commutator map factors through a well-defined map $G/Z(G) \times G/Z(G) \to [G,G]$ as stated. Furthermore, $\varpi([g],[h]) = \varpi([h],[g])^{-1}$ for all $g,h \in G$. Note that $\varpi([g],[h]) = 1$ for all $[h] \in G/Z(G)$ if and only if $[g] = 0$.
A group $G$ has nilpotency class $2$ if and only if $[G,[G,G]] = 1$ if and only if $[G,G]$ lies in $Z(G)$. This shows that $(1) \Longleftrightarrow (2)$.  Assume (2). We only need to show that $\varpi$ is bimultiplicative. Write
\begin{align*}
\varpi([g_1][g_2],[h]) = [g_1 g_2,h] & = g_1 g_2 h (g_1 g_2)^{-1} h^{-1}\\
& = g_1 g_2 h g_2^{-1} g_{1}^{-1} h^{-1}\\
& = g_1(h g_1^{-1} h^{-1})(h g_1 h^{-1})g_2hg_{2}^{-1}(h^{-1} h)g_{1}^{-1} h^{-1}\\
& = [g_1, h] (hg_1 h^{-1}) [g_2, h] (h g_{1}^{-1} h^{-1}) \\
& = [g_1, h] (hg_1 h^{-1}) [g_{2},h] (hg_1 h^{-1})^{-1} \\
& = [g_1, h] [g_2,h] \\
& = \varpi([g_1],[h])\, \varpi([g_2],[h]).
\end{align*}
Here we used the fact that $[G,G] \subseteq Z(G)$ in the penultimate equality.
Since $\varpi$ is alternating with values in $[G,G]$ which is Abelian, $\varpi$ is antisymmetric and 
\begin{equation*} \begin{split} \varpi([h],[g_1][g_2]) = \varpi([g_1][g_2],[h])^{-1} & = 
\varpi([g_1],[h])^{-1}  \varpi([g_2],[h])^{-1} \\
& = \varpi([h],[g_1])\, \varpi([h],[g_2]). \end{split}
\end{equation*} This proves (3). Conversely, assume (3). Let $g,h,k \in G$. Then rewinding the previous identity, we find that
\begin{align*}
k[g,h]k^{-1} & = [h^{-1}kh,h]^{-1}[h^{-1}khg,h] & \\
& = \varpi([h^{-1}kh],[h])^{-1}\varpi([h^{-1}khg],[h])& \\
& = \varpi([h]^{-1}[k][h],[h])^{-1}\varpi([h]^{-1}[k][h][g],[h]) & \\
& = \varpi([g],[h]) & {\hbox{(bimultiplicativity)}}\\
& = [g,h]. &
\end{align*}

\end{proof}

We state the main result of this note.



\begin{theorem} \label{th:main_th}
 Let $1 \to C \to G \to A \to 1$ be a central extension of an abelian group $A$. Assume that as ${\mathbb{Z}}$-modules, $C$ is injective or that $A$ is projective. Then

\begin{enumerate}

\item[(1)] There is a bimultiplicative map $\beta:A \times A \to C$ and an isomorphism\linebreak $G \overset{\simeq}{\to} H_{\beta}(A)$ such that the following diagram is commutative:
\begin{center}\begin{tikzcd}
1 \arrow[r] & {C} \arrow[r, hook] \arrow[d,equal] & G \arrow[r] \arrow[d, "\simeq"] & A \arrow[r] \arrow[d,equal] & 1 \\
1 \arrow[r] & {C} \arrow[r,hook] & {\mathscr{H}}_{\beta}(A) \arrow[r, "p"] & A \arrow[r] &  1.\\
\end{tikzcd} \end{center} In particular, $G$ is nilpotent of nilpotency class at most two.

\item[(2)] The associated alternating pairing $\omega_\beta:A \times A \to C$ factors through $\varpi:G/Z(G) \times G/Z(G) \to [G,G]$ by the following commutative diagram:
$$\begin{tikzcd}
{A \times A} \arrow[dr,"\omega_{\beta}"] \arrow[dd,two heads] & & \\
&  {[G,G]} \arrow[r,hook] & C \\
G/Z(G) \times G/Z(G) \arrow[ru,"\varpi"'] &  & 
\end{tikzcd}$$

\item[(3)] The following assertions are equivalent:
\begin{enumerate}
\item[(a)] The extensions $H_{\beta}(A)$ and $H_{\beta'}(A)$ are equivalent.
\item[(b)] The map $(x,y) \mapsto \beta(x,y)\beta'(x,y)^{-1}$ is symmetric.
\item[(c)] $\omega_{\beta} = \omega_{\beta'}$.
\end{enumerate}

\end{enumerate}
\end{theorem}
%
%

The proof is presented in \S \ref{sec:proof} using material from \S \ref{sec:Heisenberg_groups_as_ext}.
There are two special cases when Theorem~\ref{th:bilinear_rep} applies. The first one was proved in \cite[Th.~10.17]{CST}.

\begin{corollary_th}
Let ${\mathbb{U}}(1) = \{ z \in {\mathbb{C}} \ | \ |z| = 1 \}$.
For any central extension $1 \to {\mathbb{U}}(1) \to G \to A \to 1$ of a finite abelian group $A$, the conclusions $(1), (2)$  and $(3)$ of Theorem~$\ref{th:bilinear_rep}$ hold.
\end{corollary_th}

\begin{remark}
It is sometimes useful to replace ${\mathbb{U}}(1)$ by ${\mathbb{Q}}/{\mathbb{Z}}$ (which is also divisible).
\end{remark}

\begin{corollary_th}
For any central extension $1 \to C \to G \to A \to 1$ of a free abelian group $A$, the conclusions $(1), (2)$  and $(3)$ of Theorem~$\ref{th:bilinear_rep}$ hold. 
\end{corollary_th}

\noindent{\textbf{Acknowledgements}}. This text originated as a question about whether any Heisenberg extension up to equivalence has the form $H_{\beta}(A)$ above (where the cocycle $\beta$ is bimultiplicative) raised during a conversation with Paolo Farina. 

\section{Heisenberg groups as extensions} \label{sec:Heisenberg_groups_as_ext}

Let $A$ and $C$ be Abelian groups. We shall make the following assumption on $A$ and $C$ (as ${\mathbb{Z}}$-modules): $A$ is projective (for instance $A$ is a lattice) or $C$ is injective (for instance $C$ is ${\mathbb{C}}^{\times}$).

A \emph{$C$-valued $2$-cocycle} $c$ is a map $A \times A \to C$ such that 
\begin{equation}
c(x,1) = c(1,x) = 1 \label{eq:cocycle1}
\end{equation}
and
\begin{equation}
c(x,y)\ c(xy,z) = c(y,z) \ c(x,yz).
\label{eq:cocycle2}
\end{equation} 
For instance, a bimultiplicative map $c:A \times A \to C$ is a $C$-valued $2$-cocycle.
More examples of $C$-valued $2$-cocycles arise in the context of extensions of $A$ by $C$ which we discuss below.
A $C$-valued $2$-cocycle $c:A \times A \to C$ is {\emph{symmetric}} is $c(x,y) = c(y,x)$ for all $x, y \in A$. The {\emph{trivial}} $2$-cocycle is the map $e:A \times A \to C$ defined by $e(x,y) = 1_{C}$ for all $x,y \in A$. Clearly the product of two $2$-cocycles defined by pointwise multiplication is again a $2$-cocycle. Similarly the inverse of a $2$-cocycle $c:A \times A \to C$ is its pointwise inverse $c^{-1}$ defined by $c^{-1}(x,y) = c(x,y)^{-1}$. Associativity of product follows from associativity in $C$. We conclude that the set ${\mathscr{C}}^2(A,C)$ of $C$-valued $2$-cocycles is a group.
Let us denote

${\mathscr{C}}^2_{\rm{b}}(A,C)$ the subgroup of $C$-valued bimultiplicative maps $A \times A \to C$, 

${\mathscr{C}}^{2}_{\rm{s}}(A,C)$ the subgroup of $C$-valued symmetric $2$-cocycles,

${\mathscr{C}}^2_{\rm{sb}}(A,C)$ the subgroup of $C$-valued symmetric bimultiplicative maps $A \times A \to C$.

They a priori fit into the commutative diagram
$$
\begin{tikzcd}
& {\mathscr{C}}_{\rm{s}}^2(A,C) \arrow[rd, hookrightarrow] & \\
{\mathscr{C}}_{\rm{sb}}^2(A,C) \arrow[rd, hookrightarrow] \arrow[ru, hookrightarrow] & & {\mathscr{C}}^2(A,C) \\
& {\mathscr{C}}_{\rm{b}}^2(A,C) \arrow[ru, hookrightarrow] & 
\end{tikzcd}
$$

\begin{lemma} \label{lem:antisymmetric_bilinear}
Let $c \in {\mathscr{C}}^2(A,C)$. The map $\omega_{c}:A \times A \to C, (x,y) \mapsto c(x,y) \ c(y,x)^{-1}$ is alternating bimultiplicative.
\end{lemma}

We give a computational proof. Another proof is given after the cocycle $c$ is interpreted as the cocycle of an appropriate extension.

\begin{proof}
The map $\omega_{c}$ is obviously alternating. We shall prove that $\omega_{c}(xy,z) = \omega_{c}(x,z)\, \omega_{c}(y,z)$. The proof is similar for the other argument. We write
\begin{align*}
\omega_{c}(xy,z) & = c(xy,z)\, c(z,xy)^{-1} & \\
& = c(x,y)^{-1}\, c(y,z)\, c(x,yz)\ c(x,y)\, c(zx,y)^{-1}\, c(z,x)^{-1} & {\hbox{by}}\ (\ref{eq:cocycle2})\\
& = c(y,z)\, c(x,yz)\, c(zx,y)^{-1}\, c(z,x)^{-1} & \\
& = \underbrace{c(y,z)c(z,y)^{-1}}_{=\omega_{c}(y,z)}c(z,y) \, c(x,yz)\, c(zx,y)^{-1}\, c(x,z)^{-1} \underbrace{c(x,z)c(z,x)^{-1}}_{=\omega_{c}(x,z)} & 
\end{align*}
It remains to see that the product of the four central  terms is trivial. Since $A$ is abelian, $c(x,yz) = c(x,zy)$ so (\ref{eq:cocycle2}) applies also when $y$ and $z$ are switched:
$$ c(x,z)\, c(xz,y) = c(z,y)\, c(x,zy).$$
This is the desired result.
\end{proof}

\begin{definition}
Let $f:A \to G$ be any map from the abelian group $A$ to a group $G$. The ``morphism defect'' of $f$ is the map $\Delta {f}:A \times A \to G$ defined by $\Delta {f}(x,y) = f(xy) f(x)^{-1} f(y)^{-1}$.
A \emph{$C$-valued $2$-coboundary} $c$ is a $2$-cocycle $c:A \times A \to C$ such that there exists a map $f:A \to C$ such that
\begin{equation}
f(1) = 1\ \ \ {\rm{and}}\ \ c(x,y) = \Delta {f}(x,y)
\label{eq:cob}
\end{equation} 
\end{definition}

Let us denote by ${\mathscr{B}}^2(A,C)$ the group of $C$-valued $2$-coboundaries. It is an immediate observation that a $2$-coboundary is a symmetric $2$-cocycle: $${\mathscr{B}}^2(A,C) \subseteq {\mathscr{C}}_{\rm{s}}^2(A,C).$$

Given an abstract central extension ${\mathscr{H}}(A,C)$ of $A$ by $C$, that is, a short exact sequence $1 \to C \to {\mathscr{H}}(A,C) \overset{p}{\to} A \to 1$, one natural way to produce a $2$-cocycle $c:A \times A \to C$ is to measure the ``morphism defect'' of a section $s:A \to {\mathscr{H}}(A,C)$ of the projection map $p:{\mathscr{H}}(A,C) \to A$. Namely, given a set-theoretic section $s:A \to {\mathscr{H}}(A,C)$ such that $s(1_A) = 1$, the map defined by $$c_{s}:A \times A \to {\mathscr{H}}(A,C), \ c_{s}(x,y) = s(xy)s(x)^{-1}s(y)^{-1}$$ takes values in $C \subset {\mathscr{H}}(A,C)$ and defines a $C$-valued $2$-cocycle. 

\begin{remark}
The $C$-valued $2$-cocycle $c_{s}$ is not necessarily a $C$-valued $2$-coboundary. It is if the set-theoretic section $s$ takes values in $C$ rather than in the larger group ${\mathscr{H}}(A,C)$.
\end{remark}

%


Conversely, a $C$-valued $2$-cocycle $c:A \times A \to C$ provides a central extension ${\mathscr{H}}(A,C,c)$ fitting in the short exact sequence
\begin{equation} \begin{tikzcd}
0 \arrow[r] & C \arrow[r] & {\mathscr{H}}(A,C,c) \arrow[r] & A \arrow[r] & 0.
\end{tikzcd} \label{eq:ext} \end{equation}
The group ${\mathscr{H}}(A,C,c)$ is defined as the set $C \times A$ with group law
\begin{equation}
(t,a) \cdot (t',a') = (tt'c(a,a'), aa'),\ \ \ t,t' \in C,\ a,a' \in A. \label{eq:general_group_law}
\end{equation}
The monomorphism in the short exact sequence is the natural inclusion map \linebreak $C \to {\mathscr{H}}(A,C,c), \  t \mapsto (t,1_{A})$.
The epimorphism is the natural projection map ${\mathscr{H}}(A,C,c) \to A, \ (t,x) \mapsto x$.
It follows from $(\ref{eq:cocycle1})$ that $(1_{C}, 1_{A})$ is the unit element of ${\mathscr{H}}(A,C,c)$. 
The normalizing cocycle relation $(\ref{eq:cocycle1})$ also ensures that $C \times \{1_{A} \} \subseteq Z({\mathscr{H}}(A,C,c))$.
The cocycle relation $(\ref{eq:cocycle2})$ ensures associativity. A direct computation shows that $[{\mathscr{H}}(A,C,c), {\mathscr{H}}(A,C,c)] \subseteq C \times \{ 1_{A} \}$. Therefore (Prop.~\ref{prop:nilpotent2}) ${\mathscr{H}}(A,C,c)$ is a nilpotent group of nilpotency class (at most) $2$. 
 Furthermore, the $2$-cocycle $c$ is recovered as the ``morphism defect'' of the section $s:A \to {\mathscr{H}}(A,C,c)$ defined by $s(a) = (1_{C},a)$, $a \in A$. Indeed, $c = \Delta {s}$. Therefore, any $C$-valued $2$-cocycle is realized as the ``morphism defect'' of a section of the projection morphism ${\mathscr{H}}(A,C,c) \to A$ where ${\mathscr{H}}(A,C,c)$ is the extension associated to the $2$-cocycle $c$. 
This leads to another proof of Lemma~\ref{lem:antisymmetric_bilinear}.

\begin{proof}[Alternative proof of Lemma~\ref{lem:antisymmetric_bilinear}]
By the previous discussion, there is a section $s:A \to {\mathscr{H}}(A,C,c)$ such that $c(x,y)=s(xy)s(x)^{-1}s(y)^{-1}$, $x,y \in A$. Therefore
\begin{align*}
c(x,y)c(y,x)^{-1} 
& = s(xy)s(x)^{-1}s(y)^{-1}(s(yx)s(y)^{-1}s(x)^{-1})^{-1} \\
& = s(xy)s(x)^{-1}s(y)^{-1}s(x)s(y)s(yx)^{-1}  \\
& = s(xy)[s(x)^{-1},s(y)^{-1}]s(yx)^{-1}  \\
& = s(xy) [s(x)^{-1},s(y)^{-1}] s(xy)^{-1}  \\
& = [s(x)^{-1},s(y)^{-1}]
\end{align*}
Here we used that $A$ is abelian in the penultimate equality and that $[s(x)^{-1}, s(y)^{-1}]$ lies in the central (normal) subgroup $C$. In terms of Proposition~\ref{prop:nilpotent2}, we have proved that 
\begin{equation}
\omega_{c}(x,y) = \varpi(s(x)^{-1},s(y)^{-1}) = \varpi(s(x),s(y)). \label{eq:rels_omegas}
\end{equation} 
(The last equality in virtue of $\varphi$ being bimultiplicative again by Prop.~\ref{prop:nilpotent2}.)
Next,
\begin{align*}
 \omega_{c}(xx',y) = \left[(c(x,x')s(x)s(x'))^{-1},s(y)^{-1}\right] & = [s(x')^{-1}s(x)^{-1}c(x,x')^{-1}, s(y)^{-1}] \\
& = [s(x')^{-1}s(x)^{-1},s(y)^{-1}]
\end{align*} where for the last equality we used the fact that $c(x,x')^{-1} \in C \subseteq Z({\mathscr{H}}(A,C,c))$. It follows from Prop.~\ref{prop:nilpotent2} that
\begin{align*}
 \omega_{c}(xx',y) = \varpi(s(x')^{-1}s(x)^{-1},s(y)^{-1}) & = 
\varpi(s(x')^{-1}, s(y)^{-1}) \, \varpi(s(x)^{-1},s(y)^{-1}) \\
& = 
\omega_{c}(x',y) \, \omega_{c}(x,y) \\
& = \omega_{c}(x,y)\, \omega_{c}(x',y).
\end{align*}
This is the desired result.
\end{proof}

%
%
%

\begin{remark} \label{rem:basic_rem}
The constructions recalled above, of a $2$-cocycle $c_{s}$ from an abstract central extension ${\mathscr{H}}(A,C)$ and a section $s:A \to {\mathscr{H}}(A,C)$ on the one hand, and of the central extension ${\mathscr{H}}(A,C,c)$ from a $2$-cocycle $c$ on the other hand, are inverse of each other in the sense that the extensions ${\mathscr{H}}(A,C)$ and ${\mathscr{H}}(A,C,c_{s})$ are equivalent, that is, fit into the commutative diagram with exact rows
\begin{equation} \begin{tikzcd}
0 \arrow[r] & C \arrow[r] \arrow[d, equal ] & {\mathscr{H}}(A,C) \arrow[r] \arrow[d,"\simeq"] & A \arrow[r] \arrow[d, equal] & 0 \\
0 \arrow[r] & C \arrow[r] & {\mathscr{H}}(A,C,c_s) \arrow[r] & A \arrow[r] & 0.
\end{tikzcd} \label{eq:equiv_ext} 
\end{equation}
An explicit isomorphism is given by the map ${\mathscr{H}}(A,C,c_{s}) \to {\mathscr{H}}(A,C), \ (t,a) \mapsto ts(a).$
\end{remark}

A consequence of this and the discussion above is another characterization of nilpotent groups of nilpotency class (at most) $2$:

\begin{prop}
A group $G$ is nilpotent of class (at most) $2$ if and only if it sits in a central extension $1 \to C \to G \to A \to 1$ where $A$ is an abelian group.
\end{prop}

\begin{definition}
The {\emph{second cohomology group}} $H^{2}(A,C)$ is defined as $H^{2}(A,C) = {\mathscr{C}}^2(A,C)/{\mathscr{B}}^2(A,C)$.
\end{definition}

\begin{remark}
The property that a (necessarily symmetric) bimultiplicative pairing $c:A \times A \to C$ is a $2$-coboundary is equivalently expressed by the property that $c:A \times A \to C$ has a quadratic refinement in the sense of \cite{DM}. In particular, this property holds if $H^{2}(A,C) = 0$.
\end{remark}

It is a fundamental fact of homological algebra that the set of equivalence classes of extensions (\ref{eq:ext}) is in bijective correspondence with $H^{2}(A,C)$ (see e.g., \cite[Chap.~5]{Rotman}). Thus two extensions ${\mathscr{H}}(A,C,c)$ and ${\mathscr{H}}(A,C,c')$ are equivalent (as in the sense defined above)
if and only if $c\, c'^{-1} \in {\mathscr{B}}^2(A,C)$.

We state the main result of this section.

\begin{theorem} \label{th:bilinear_rep}
Every $2$-cohomological class $[c] \in H^2(A,C)$ has a representative in ${\mathscr{C}}^2_{\rm{b}}(A,C)$. 
Two $2$-cocycles $c, c' \in {\mathscr{C}}^{2}(A,C)$ are cohomologous if and only if $\omega_{c} = \omega_{c'}$.
\end{theorem} 

The particular case when $A$ is finite abelian and $C = {\mathbb{C}}$ is proved in \cite[\S 10.1]{CST}.
The proof of Theorem~\ref{th:bilinear_rep} is based on the following observation.

\begin{prop}\label{prop:symmetric_implies_coboundary}
A $2$-cocycle in ${\mathscr{C}}(A,C)$ is symmetric if and only if it is a $2$-coboundary: ${\mathscr{C}}_{s}^{2}(A,C) = {\mathscr{B}}^{2}(A,C)$
\end{prop}

\begin{remark} In this remark we use additive notation.
It follows from Prop.~\ref{prop:symmetric_implies_coboundary} that a symmetric bilinear pairing $c:A \times A \to C$ has a quadratic refinement $q:A \to C$ such that $c = \Delta q$. This fact is nontrivial. For instance, let $A = C = {\mathbb{Z}}$ and $c(x,y) = xy$. Clearly there is no \emph{homogeneous} quadratic form $q$ over ${\mathbb{Z}}$ such that $c = \Delta q$. However, as implied by Prop.~\ref{prop:symmetric_implies_coboundary}, there is a nonhomogenous quadratic map $q:{\mathbb{Z}} \to {\mathbb{Z}}$ such that $c = \Delta q$. In additive notation, one verifies that $q(x) = \frac{x^2}{2} - \frac{x}{2}$ defines a quadratic map $q:{\mathbb{Z}} \to {\mathbb{Z}}$ and $\Delta q = c$.
\end{remark}

Since Prop.~\ref{prop:symmetric_implies_coboundary} expresses a basic fact, we shall give two proofs. The first proof relies on some basic commutative algebra and seems new. The second proof, valid only if $C = {\mathbb{F}}^{\times}$ where ${\mathbb{F}}$ is an algebraically closed field, relies on representation theory and is due to \cite{CST}.

\begin{proof}[First proof of Prop.~$\ref{prop:symmetric_implies_coboundary}$] 
A $2$-coboundary is symmetric. Let us prove the converse.
Let $c \in {\mathscr{C}}_{s}^{2}(A,C)$. Consider the corresponding Heisenberg extension ${\mathscr{H}}(A,C,c)$. Observe that since $c$ is symmetric (and $A$ abelian), the group law (\ref{eq:general_group_law}) on ${\mathscr{H}}(A,C,c)$ is commutative. So we have a short exact sequence of abelian groups
$$ 1 \to C \to {\mathscr{H}}(A,C) \overset{p}{\to} A \to 1 $$
Since $C$ is injective or $A$ is projective, the short exact sequence splits (as a sequence of ${\mathbb{Z}}$-modules): there is a homomorphism $\sigma:A \to {\mathscr{H}}(A,C)$ such that $p \circ \sigma = {\rm{Id}}_{A}$. This last property implies that there is a map $q_{\sigma}:A \to C$ such that $\sigma(x) = (q_{\sigma}(x),x)$, $x \in A$ satisfying
$$ (q_{\sigma}(xy), xy) = \sigma(xy) = \sigma(x)\sigma(y) = (q_{\sigma}(x),x) \cdot (q_{\sigma}(y),y) = (q_{\sigma}(x)q_{\sigma}(y)c(x,y), xy),$$
hence $q_{\sigma}(xy) = q_{\sigma}(x)q_{\sigma}(y)c(x,y)$, i.e., $c \in {\mathscr{B}}^{2}(A,C)$.
\end{proof}

For the second proof, we need to define a projective representation.


\begin{definition}
A projective representation associated to a $2$-cocycle $c:A \times A \to C$ is a map $\rho_{c}:A \to {\rm{GL}}(V)$ such that $\rho_{c}(xy) = c(x,y)^{-1}\, \rho_{c}(x) \, \rho_{c}(y)$, $x, y \in A$. 
\end{definition}

The reason for the convention in the definition should be clear after the following proposition.

\begin{prop} \label{prop:reps}
A projective representation $\rho:A \to {\rm{GL}}(V)$ with $2$-cocycle $c:A \times A \to C$ gives rise to a linear representation $\tilde{\rho}:{\mathscr{H}}(A,C,c) \to {\rm{GL}}(V)$ by
$$ \tilde{\rho}(t,x) = t\, \rho(x), \ t \in C, \ x \in A.$$
Conversely, a linear representation $\sigma:{\mathscr{H}}(A,C,c) \to {\rm{GL}}(V)$ such that $\sigma(t,1) = t\, {\rm{Id}}_{V}$, $t \in C$, gives rise to a projective representation $\sigma':A \to {\rm{GL}}(V)$ with $2$-cocycle $c$ by $\sigma'(x) = \sigma(1,x)$, $x \in A$.
The map $\rho \mapsto \tilde{\rho}$ is a bijection with inverse $\sigma \mapsto \sigma|_{1 \times A}$ between the set of all projective 
representations of $A$ with cocycle $c$ and the set of all linear representations $\theta$ of ${\mathscr{H}}(A,C,c)$ such that $\theta(t,1) = t\, {\rm{Id}}_{V}$ for all $t \in C$. Furthermore, the map and its inverse preserve unitarity, irreducibility and equivalence.
\end{prop}

\begin{proof}[Proof of Prop. $\ref{prop:reps}$]
For the first statement,
\begin{align*}
 \tilde{\rho}((t,x)(t',y)) = \tilde{\rho}(tt'c(x,y),xy) = tt'c(x,y) \rho(xy)  & = tt'c(x,y)\, c(x,y)^{-1} \rho(x) \rho(y) \\
& = t\rho(x) \, t'\rho(y) = \tilde{\rho}(t,x)\, \tilde{\rho}(t',y).
\end{align*}
Clearly, $\tilde{\rho}(t,1) = t \, \rho(1) = t\, {\rm{Id}}_{V}$.
Conversely, 
\begin{align*}
\sigma'(x)\sigma'(y) = {\sigma}(1,x){\sigma}(1,y) = \sigma((1,x)(1,y)) & = 
\sigma(c(x,y),xy) \\
& = \sigma((c(x,y),1)(1,xy)) \\
& = \sigma(c(x,y),1) \ \sigma(1,xy)\\
& = c(x,y)\, {\rm{Id}}_V \ \sigma'(xy) \\
& = c(x,y) \ \sigma'(xy).
\end{align*}
So $\sigma'$ is a projective representation associated to the cocycle $c$.
The remaining statements are clear.
\end{proof}

\begin{remark}
Now it is clear that our definition of a projective representation associated to a $2$-cocycle $c$ fits our definition of the law group of the extension ${\mathscr{H}}(A,C,c)$.
\end{remark}

\begin{proof}[Second proof of Prop.~$\ref{prop:symmetric_implies_coboundary}$ {\cite[Theorem~10.17]{CST}}] 
Let $c \in {\mathscr{C}}_{s}^{2}(A,C)$ a $2$-cocycle giving rise to a central extension ${\mathscr{H}}(A,C,c)$. As $C = {\mathbb{F}}^{\times} = {\rm{GL}}({\mathbb{F}})$, the identity map $\chi: C \to C$ is a one-dimensional (irreducible) representation of $C$. Regard $C$ as a subgroup of ${\mathscr{H}}(A,C,c)$. There is an induced linear representation $\pi = {\rm{Ind}}_{C}^{\mathscr{H}(A,C,c)}(\chi): {\mathscr{H}}(A,C,c) \to {\rm{GL}}(V)$ for some vector space $V$ over ${\mathbb{F}}$, such that 
$$ \pi(t,1) = t \, {\rm{Id}}_{V}, \ t \in C.$$
Since $\chi$ is irreducible on a central subgroup, the induced representation $\pi$ is irreducible.
By Prop.~\ref{prop:reps}, there exists an irreductible projective representation $\rho:A \to {\rm{GL}}(V)$ over a vector space $V$ over ${\mathbb{F}}$ with cocycle $c$:
\begin{equation}
 \rho(x)\, \rho(y) = c(x,y)\, \rho(x\, y), \ \ \ x,y \in A. \label{eq:proj_rep_with_cocycle}
 \end{equation}
Now since $c$ is symmetric, $c(x,y) = c(y,x)$ and $A$ is abelian,
$$ \rho(y)\, \rho(x) = c(y,x)\, \rho(y\, x) = c(x,y)\, \rho(y\, x) 
= c(x,y) \, \rho(x\, y) = \rho(x) \, \rho(y).$$
Thus $\rho(x)$ commutes with all $\rho(y)$, $y \in A$. Schur's lemma over ${\mathbb{F}}$ (here we use the fact that ${\mathbb{F}}$ is algebraically closed) implies that $\rho(x)$ is a homothety: there exists $\lambda(x) \in {\mathbb{F}}^{\times} (= C)$ such that $\rho(x) = \lambda(x)\, {\rm{Id}}_{V}$. But then (\ref{eq:proj_rep_with_cocycle}) implies that
$$ \lambda(x)\, \lambda(y) = c(x,y) \, \lambda(x\, y), \ \ \ x,y \in A $$
which means that $c \in {\mathscr{B}}^{2}(A,C)$.
\end{proof}


\begin{definition}
The set of alternating bimultiplicative pairings $A \times A \to C$ is denoted ${\mathscr{C}}^{2}_{\rm{altb}}(A,C)$. 
\end{definition}

\begin{remark}
Note that ${\mathscr{C}}^{2}_{\rm{altb}}(A,C)$ is a subgroup of ${\mathscr{C}}^{2}_{\rm{b}}(A,C)$. 
\end{remark}

We are now ready for the proof of Theorem~\ref{th:bilinear_rep}.

\begin{proof}[Proof of Theorem~\ref{th:bilinear_rep}]
According to Lemma~\ref{lem:antisymmetric_bilinear}, there is a well-defined map $${\mathscr{C}}^{2}(A,C) \to {\mathscr{C}}^{2}_{\rm{altb}}(A,C), \ c \mapsto \omega_{c}.$$
Using commutativity of $A$, it is readily verified that this map is a group homomorphism. Furthermore, its kernel is clearly ${\mathscr{C}}^{2}_{\rm{s}}(A,C)$. We claim that the restriction of this map on ${\mathscr{C}}^{2}_{\rm{b}}(A,C)$ is already onto. Indeed, let us consider an alternate bimultiplicative pairing $g:A \times A \to C$. Choose a minimal system of symplectic generators $a_1, \ldots, e_{2n} \in A$ such that $g(a_{i}, a_{i}) = 0$, $i=1, \ldots, 2n$. Define a bimultiplicative map $c:A \times A \to C$ by $$
\left\{ \begin{array}{ll}
c(a_{i}, a_{j}) = g(a_{i}, a_{j}) & 1 \leq i < j \leq 2n \\
c(a_{j},a_{i}) = 1 &  1 \leq i \leq j \leq 2n.
\end{array} \right.$$
Then $\omega_{c} = g$. We therefore have the following commutative diagram
$$ \begin{tikzcd}
& {\mathscr{C}}^{2}_{\rm{sb}}(A,C) \arrow[r,hook] 
\arrow[d,hook] 
& {\mathscr{C}}^{2} _{\rm{b}}(A,C) \arrow[r,"\omega_{\bullet}|"] \arrow[d,hook] & {\mathscr{C}}^{2}_{\rm{altb}}(A,C) \arrow[r] \arrow[d,equal] & 1 \\
{\mathscr{B}}^{2}(A,C) \arrow[r,equal] & {\mathscr{C}}^{2}_{\rm{s}}(A,C) \arrow[r,hook] & {\mathscr{C}}^{2}(A,C) \arrow[r,"\omega_{\bullet}"] & {\mathscr{C}}^{2}_{\rm{altb}}(A,C) \arrow[r] & 1.
\end{tikzcd} $$
It follows that
$$ H^{2}(A,C) = {\frac{{\mathscr{C}}^{2}(A,C)}{{\mathscr{B}}^2(A,C)}} = {\frac{{\mathscr{C}}^{2}(A,C)}{{\mathscr{C}}^{2}_{\rm{s}}(A,C)}} = \frac{\mathscr{C}^{2}_{\rm{b}}(A,C)}{\mathscr{C}^{2}_{\rm{sb}}(A,C)} \simeq {\mathscr{C}}^{2}_{\rm{altb}}(A,C).$$
\end{proof}

\section{Corollaries}

Theorem~\ref{th:bilinear_rep} has a number of important consequences.

\begin{prop} \label{prop:2-cohomology_equals_symplectic}
The map $\omega_{\bullet}:{\mathscr{C}}^{2}(A,C) \to {\mathscr{C}}^{2}_{\rm{altb}}(A,C)$ induces an isomorphism $H^{2}(A,C) \overset{\simeq}{\to} {\mathscr{C}}^{2}_{\rm{altb}}(A,C)$.
\end{prop}

\begin{prop} \label{prop:general_heisenberg_has_bimultiplicative_pairing}
Any abstract Heisenberg extension ${\mathscr{H}}(A,C)$ is equivalent to a Heisenberg group ${\mathscr{H}}_{\beta}(A)$ for some bimultiplicative pairing $\beta:A \times A \to C$. Furthermore, any two Heisenberg extensions $H_{\beta}(A)$ and $H_{\beta'}(A)$ are equivalent if and only if $\omega_{\beta} = \omega_{\beta'}$ if and only if $\beta \, \beta'^{-1}$ is symmetric (bimultiplicative).
\end{prop}

\begin{proof}
Direct consequence of Theorem~\ref{th:bilinear_rep} and Prop.~\ref{prop:symmetric_implies_coboundary}. 
\end{proof}

\begin{prop}\label{prop:cyclic_trivial_2-cohomology}
If $A$ is cyclic then $H^{2}(A,C) = 0$.
\end{prop}

\begin{proof}
The hypothesis implies that there is no nontrivial alternating bimultiplicative pairing on $A$, so ${\mathscr{C}}^{2}_{\rm{altb}}(A,C)$ is trivial, hence, by Prop.~\ref{prop:2-cohomology_equals_symplectic} above, the result.
\end{proof}

\section{Proof of the main theorem} \label{sec:proof}


According to Remark~\ref{rem:basic_rem}, there exists a section $s:A \to G$ and a $C$-valued $2$-cocycle $c = c_{s}:A \times A \to C$ such that the diagram
$$\begin{tikzcd}
1 \arrow[r] & C \arrow[r, hook] \arrow[d,equal] & G \arrow[r] \arrow[d, "\simeq"] & A \arrow[r] \arrow[l, bend right, "s"'] \arrow[d,equal] & 1 \\
1 \arrow[r] & C \arrow[r,hook] & {\mathscr{H}}(A,C,c) \arrow[r, "p"] & A \arrow[r] &  1
\end{tikzcd}$$
is commutative, i.e. $G$ and ${\mathscr{H}}(A,C,c)$ are equivalent as extensions. By Prop.~\ref{prop:general_heisenberg_has_bimultiplicative_pairing}, ${\mathscr{H}}(A,C,c)$ is equivalent to $H_{\beta}(A)$ for some bimultiplicative pairing $\beta:A \times A \to C$. Therefore we can complete the commutative diagram of extension equivalences
$$\begin{tikzcd}
1 \arrow[r] & C \arrow[r, hook] \arrow[d,equal] & G \arrow[r] \arrow[d, "\simeq"] & A \arrow[r] \arrow[d,equal] & 1 \\
1 \arrow[r] & C \arrow[r,hook] \arrow[d,equal] & {\mathscr{H}}(A,C,c) \arrow[r, "p"] \arrow[d,"\simeq"] & A \arrow[r] \arrow[d,equal] &  1\\
1 \arrow[r] & C \arrow[r,hook] & {\mathscr{H}}_{\beta}(A) \arrow[r, "p"] & A \arrow[r] &  1
\end{tikzcd}$$
 This gives the first commutative diagram of the Theorem. Then
\begin{align*}
\omega_{\beta}(x,y) & = \omega_{c}(x,y) & \hbox{by Theorem~\ref{th:bilinear_rep}}\\ 
& = \varpi(s(x),s(y)) & \hbox{by $(\ref{eq:rels_omegas})$} \\
& = \varpi(ts(x),t's(y)) & \hbox{by Prop.~$\ref{prop:nilpotent2}$} \\
& = \varpi(\varphi(t,x),\varphi(t',y)),
\end{align*}
where $\varphi$ is the isomorphism ${\mathscr{H}}(A,C,c) \to G, 
(t,x) \mapsto ts(x)$ for some section $s$ such that $c = c_s$ (cf. Remark~\ref{rem:basic_rem}). This shows that the second diagram is commutative as well and proves the first statement. The last statement is a restatement of Proposition~\ref{prop:general_heisenberg_has_bimultiplicative_pairing} above.

\bibliographystyle{amsalpha}

\end{document}